\documentclass[reqno,12pt]{article}
\usepackage{a4wide}
\usepackage{amsmath}
\usepackage{amssymb}
\usepackage{amsthm}
\usepackage{mathrsfs}
\usepackage[utf8]{inputenc}
\usepackage{graphicx}

\numberwithin{equation}{section}

\newtheorem{theo}{Theorem}

\newtheorem{coro}{Corollary}

\newtheorem{defi}{Definition}

\theoremstyle{remark}
\newtheorem*{Remark}{Remark}

\def\ch{\chi}

\def\({\left(}
\def\){\right)}
\def\[{\left[}
\def\]{\right]}

\def\Qb{\overline{\mathbb Q}}

\newcommand{\N}{\mathbb{N}}
\newcommand{\Z}{\mathbb{Z}}
\newcommand{\Q}{\mathbb{Q}}
\newcommand{\R}{\mathbb{R}}
\newcommand{\C}{\mathbb{C}}

\newcommand{\Qbar}{\overline{\mathbb Q}}

\newcommand{\ord}{{\rm ord}}

\newcommand{\eps}{\varepsilon}

\newcommand{\moins}{\setminus}

\newcommand{\calD}{{\cal {D}}}

\newcommand{\gdo}{\mathcal{O}}

\newcommand{\pe}[1]{\lfloor #1 \rfloor}
\renewcommand{\ch}{c_1}
\newcommand{\cs}{c_5}
\newcommand{\cn}{c_6}
\newcommand{\cdix}{c_{7}}

\newcommand{\hz}{h_0}
\newcommand{\hun}{h_1}
\newcommand{\hd}{h_2}

\newcommand{\bz}{b}

\newcommand{\notn}{{\mathcal N}_b(\xi,t,n)}
\newcommand{\notnun}{{\mathcal N}_b(\xi,1,n)}
\newcommand{\notnbis}{{\mathcal N}_{\bz}\big(F(a/\bz^s),t,n\big)}

\begin{document}

\title{Rational approximation to values of $G$-functions, and their expansions in integer bases}

\author{S. Fischler and T. Rivoal}
\date{\today}

\maketitle

\begin{abstract} We prove a general effective result concerning   approximation  of irrational values at rational points $a/b$ of any $G$-function $F$ with rational Taylor coefficients by fractions of the form $n/(B\cdot b^m)$, where the integer $B$ is fixed. 
As a corollary, we show that if $F$ is not in $\mathbb Q(z)$, for any $\eps>0$ and any $b$ and $m$ large enough with respect to $a$,  $\eps$ and $F$, then  $|F(a/b)-n/b^m|\ge 1/b^{m(1+\eps)}$ and $F(a/b)\notin\mathbb Q$.  
This enables us to obtain a new and effective result on  repetition of patterns in the $b$-ary expansion of $F(a/b^s)$ for any $b\geq 2$. In particular, defining $\mathcal{N}(n)$ as the number of consecutive equal digits in the  $b$-ary  expansion of $F(a/b^s)$ starting from the $n$-th digit, we prove that $\limsup_n \mathcal{N}(n)/n\le \eps$ provided the integer $s\ge1$ is such that $b^s$ is  large enough with respect to $a$,  $\eps > 0$ and $F$. This improves over the previous bound $1+\eps$, that can be deduced from the work of  Zudilin.

Our crucial ingredient is the  use of {\em non-diagonal} simultaneous Pad\'e type approximants  for any given family of $G$-functions solution of a differential system, in a construction \`a la Chudnovsky-Andr\'e. This idea was introduced by  Beukers in the particular case of the function $(1-z)^\alpha$ in his study of the generalized Ramanujan-Nagell equation, and we use it in its full generality here.
In contrast with the classical Diophantine ``competition'' between $E$-functions and $G$-functions,  similar results are  still not known for a single transcendental value of an  $E$-function  at  a rational point, not even for the exponential function.
\end{abstract}

\section{Introduction}
This paper deals with  approximations of 
values of $G$-functions at rational points by rational numbers with denominator a power 
of a fixed integer; an important motivation is that periods are conjecturally values of $G$-functions (see \cite[Section 2.2]{firi}). Before stating our results, we recall some important results in the Diophantine 
theory of $G$-functions, as well as of $E$-functions, even though no new result will be given for the latter. Throughout the paper we fix an embedding of 
$\Qbar$ into $\C$. 
\begin{defi} \label{def:gfunc}
A $G$-function $F$ is a power series $F(z)=\sum_{n=0}^{\infty} a_n z^n$ 
such that the coefficients $a_n$ are algebraic numbers and there exists $C>0$ such that, for any $n\geq 1$:

$(i)$ the maximum of the moduli of the conjugates of $a_n$ is $\leq C^n$.

$(ii)$ there exists a sequence of rational integers $d_n$, with $\vert d_n \vert \leq C^n$,  such that 
$d_na_m$ is an algebraic integer for all~$m\le n$.

$(iii)$ $F(z)$ satisfies a  homogeneous linear differential equation with 
coefficients in $\Qb(z)$.

\noindent An $E$-function is a power series 
$F(z)=\sum_{n=0}^{\infty} \frac{a_n}{n!} z^n$ such that $\sum_{n=0}^{\infty} a_n z^n$ is a $G$-function.
\end{defi}
Siegel's original definition~\cite{siegel} of $E$ and $G$-functions is slightly more general but it is believed to define the same functions as above. It is a fact that the Diophantine theory of $G$-functions is not as fully 
developped as that of $E$-functions. There is no general theorem about 
the transcendence of values of $G$-functions, but results like the following one, due to Chudnovsky~\cite{chud1}. 

\medskip

{\em Let $N\ge2$ and $Y(z)={}^t(F_1(z), \ldots, F_N(z))$ be a vector of $G$-functions solution of 
a differential system $Y'(z)=A(z)Y(z)$, where $A(z)\in M_N(\Qb(z))$. Assume that $F_1(z)$, \ldots, $F_N(z)$ are $\C(z)$-algebraically independent. 
Then for any $d$, there exists $C=C(Y,d)>0$ such that, for any algebraic number $\alpha\neq 0$ of degree $d$ with  $\vert \alpha \vert <\exp(-C\log\left(H(\alpha)\right)^{\frac{4N}{4N+1}}),$ there does not exist a polynomial relation between the  values $F_1(\alpha), \ldots, F_N(\alpha)$ over $\mathbb Q(\alpha)$ of degree $d$.}

\medskip

Here, $H(\alpha)$ is the naive height of $\alpha$, i.e. the maximum of the modulus of the coefficients of the (normalized) minimal polynomial of $\alpha$ over $\mathbb Q$. Chudnovsky's theorem refines the works 
of Bombieri~\cite{bombieri} and Galochkin~\cite{galoshkin}.  Andr\'e generalized Chudnovsky's theorem to the case of an inhomogenous system $Y'(z)=A(z)Y(z)+B(z)$, 
$A(z), B(z)\in M_N(\Qb(z))$, with a similar condition on $\alpha$ and $H(\alpha)$; see~\cite[pp. 130--138]{andre} when the place $v$ is archimedean. Andr\'e and Chudnovsky's theorems are still essentially the best known today in this generality but they are far from being transcendence or algebraic independence statements. We recall that, in fact, it is not even known if there exist three algebraically independent $G$-values.~(\footnote{There exist examples of 
two algebraically independent $G$-values, for instance $\pi$ and $\Gamma(1/3)^3$, 
or $\pi$ and $\Gamma(1/4)^4$. This was first proved by Chudnovsky with a method not related to $G$-functions, but 
Andr\'e~\cite{andrecrelle} obtained a proof with certain Gauss' 
hypergeometric functions, which are $G$-functions. Andr\'e's method is very specific and has not been  
generalized.}) 
On the other hand,  
the situation is best possible for $E$-functions, by the Siegel-Shidlovsky Theorem \cite{siegel,Shidlovski}: 
except maybe for a finite set included in the set of 
singularities of a given differential system satisfied by $E$-functions, the numerical transcendence degree over $\Qbar$ of the values of the latter at a  non-zero algebraic point  is equal to their functional transcendence degree over $\Qbar(z)$.  
Beukers~\cite{beukers} was even able to describe very precisely the nature of the numerical algebraic relations when the transcendence degree is not maximal.

\medskip

A lot of work has been devoted to improvements of Chudnovsky's theorem, 
or alike, for classical  $G$-functions like 
the polylogarithms $\sum_{n=1}^\infty z^n/n^s$, or to determine weaker conditions for the irrationality of 
the values of $G$-functions at rational points. From a qualitative point of view, the result is the following.

\medskip

{\em Let $F$ be a  $G$-function with rational Taylor coefficients such that $F(z)\not\in\Q(z)$. Then there exist   positive constants $C_1$ and $C_2$, depending only on $F$, with the following property. Let $a\neq 0$ and $b \geq 1$ be integers such that 
\begin{equation} \label{eqirratf}
b > (C_1 |a|)^{C_2}  .
\end{equation}
Then $F(a/b)$ is irrational.}

\medskip

This result follows from Theorem I in \cite{chud1, chud2}, together with an irrationality measure; see also  \cite{galoshkin}. 
This measure and the value of $C_2$ have been improved by Zudilin  \cite{zudilin}, under further 
assumptions on $F$. He obtains  the following result (in a more precise form).
\begin{theo}[Zudilin~\cite{zudilin}]\label{theo:zud} Let $N\ge2$ and $Y(z)={}^t(F_1(z), \ldots, F_N(z))$ be a vector of $G$-functions solution of 
a differential system $Y'(z)=A(z)Y(z)+B(z)$, where $A(z), B(z)\in M_N(\mathbb C(z))$. Assume either  that 
$N=2$ and $1,F_1(z), F_2(z)$ are $\mathbb C(z)$-linearly independent, or that $N\ge 3$ and 
$F_1(z), \ldots ,F_N(z)$ are $\mathbb C(z)$-algebraically independent. Let $\eps > 0$,  $a\in \Z$, $a\neq 0$. Let $b$ and $q$ be 
sufficiently large positive integers, in terms of the $F_j$'s, $a$ and $\eps$; then $F_j(a/b)$ is an irrational number and for any integer $p$, we have
\begin{equation}\label{eq:zud}
\left\vert F_j \Big(\frac ab\Big) - \frac{p}{q}\right\vert \ge   \frac{1}{q^{2+\eps}}, \quad j=1, \ldots, N.
\end{equation}
\end{theo}
Zudilin's proof follows Shidlovsky's ineffective approach to zero estimates (see  \cite[p.~93, Lemma 8]{Shidlovski}). It is likely that using an effective method instead (see \cite[Appendix of Chapter~III]{andre}, \cite{BB} and \cite{ChudHouches}), one would make Theorem \ref{theo:zud} effective. 
We mention that Zudilin~\cite{zudilin2} also obtained similar irrationality measures for the values of $E$-functions at any non-zero rational point. 

\medskip

We now come to our main result. Roughly speaking, it is an improvement of Zudilin's exponent 
$2+\eps$ in~\eqref{eq:zud} when $q$ is restricted to integers of the form $b^m$.  
In this case, the exponent drops from $2+\eps$ to $1+\eps$; see Corollary~\ref{coro} with $B=1$. We first state a more precise and general version, without $\eps$, which contains an irrationality measure in disguise (see the comments following the theorem). 
\begin{theo} \label{th1}
 Let $F$ be a  $G$-function with rational Taylor coefficients and with $F(z)\not\in\Q(z)$, and $t\geq0$.
Then there exist some positive effectively computable constants $c_1$, $c_2$, $c_3$, $c_4$, depending only on $F$ (and $t$ as well for $c_3$), such that the following property holds. 
Let $a\neq 0$ and $b,B\geq 1$ be integers such that 
\begin{equation}
b > (c_1 |a|)^{c_2} \mbox{ and } B \leq b^t. \label{eqminob}
\end{equation}
Then $F(a/b)\not\in\Q$, and for any $n\in\Z$ and any $m\geq c_3 \frac{\log(b)}{\log(|a|+1)}$ we have
\begin{equation}\label{eq:mesres}
\left\vert F \Big(\frac ab\Big) - \frac{n}{B\cdot b^m}\right\vert \ge \frac{1}{B\cdot b^m \cdot (|a|+1)^{c_4m} }.
\end{equation}
\end{theo}
In the case of the dilogarithm $\textup{Li}_2(z)=\sum_{n=1}^\infty \frac{z^n}{n^2}$, our proof provides $c_1=4 e^{66}$, $c_2=12$ 
and $c_4=10^{6}$. We did not try to compute $c_3$ because it is useless for the application stated in Theorem~\ref{th2} below, but this could be done in principle. Needless to say, these values are far from   best possible but 
this is not the point of this paper. 
For related  
results, but restricted only to the $G$-functions $(1-z)^\alpha$ and $\log(1-z)$, see~\cite{BauerBennett,BugeaudBennett,beukers1} and~\cite{riv1} respectively. 
We point out that Theorem \ref{th1} is effective, because an effective zero estimate (due to Andr\'e) is used. In contrast with 
Zudilin's theorem, we only need to assume that $F(z)\not\in\mathbb Q(z)$. 

\medskip

The lower bound \eqref{eq:mesres} implies an effective irrationality measure of $F(a/b)$. Indeed, let $A  $ and $B\ge 1$ be any integers,   $t=\frac{\log(B)}{\log(b)}$ and  $m=\pe{c_3 \frac{\log(b)}{\log(|a|+1)}}+1$. The proof of Theorem~\ref{th1} shows that one may take   $c_3 = \frac43 t $ if $t$ (i.e. $B$ with our choice here) is large enough in terms of $F$.  Then, with $n=A\cdot b^m$,  Eq. 
\eqref{eq:mesres} implies that, provided $b > (c_1 |a|)^{c_2}$, we have 
\begin{equation}\label{eq:irratalt}
\left\vert F \Big(\frac ab\Big) - \frac{A}{B}\right\vert \ge \frac{\kappa}{B^{\mu}}
\end{equation}
for some  constants $\kappa,\mu>0$ that depend effectively on $a$, $b$, and $F$. The constant $\mu$ is worse that Zudilin's, at least when $b$ is large with respect to $a$, but~\eqref{eq:irratalt} applies to a larger class of $G$-functions. On the other hand, when $F(z)\not\in \mathbb Q(z)$, we can compare~\eqref{eq:irratalt} with  Chudnovsky's irrationality measure~\cite{chud1, chud2} under   assumption~\eqref{eqirratf}: our constant $c_2=3(N+2)$ in \eqref{eqminob} is slightly worse than his $C_2=N(N+1)/\eps+N+1$ (when his $\eps>0$ is large) but we have not been able to compute his value $C_1$, which depends on $F$ and $\eps$. In any case, it is difficult to compare such results in the literature as they apply to more or less  $G$-functions, to more or less values $a/b$, and give more or less close to optimal irrationality measures.

\medskip

We now extract a lower bound similar to Zudilin's measure~\eqref{eq:zud}. 
Given $\eps>0$ and assuming that  
$m \geq  2t / \eps $ and 
$b > (|a|+1)^{2 c_4/\eps }$, we derive the following corollary from Theorem~\ref{th1}. It is also effective and it will be used to prove Theorem~\ref{th2} below.
\begin{coro} \label{coro}
 Let $F$ be a   $G$-function with rational Taylor coefficients and with $F(z)\not\in\Q(z)$, $\eps > 0$,  $t\geq0$ 
and $a\in \Z$, $a\neq 0$. 
Let $b$ and $m$ be positive integers, sufficiently large in terms of $F$, $\eps$, 
  $a$   (and $t$ for $m$).  Then $F(a/b)\not\in\Q$ and for any 
integers $n$ and $B$ with $1\leq B \leq b^t$, 
we have
$$
\left\vert F \Big(\frac ab\Big) - \frac{n}{B\cdot b^m}\right\vert \ge  \frac1{b^{m(1+\eps)}}.
$$
\end{coro}
We don't know 
if some analogues of Theorem~\ref{th1} and Corollary~\ref{coro} hold when $F$ is supposed to 
be an $E$-function. This is surprising  because, as we indicated above, the Diophantine theory of 
$E$-values is much more advanced than that of $G$-values. Our method is inoperant for $E$-functions 
and we could not find any way to fix it. We explain the reason for this unusual advantage of $G$-functions in the final 
Section~\ref{sec:conclusion}. We also explain there that an analogue of Theorem~\ref{th1} holds 
for $1/F(z)$ instead of $F(z)$ under a less general assumption on the $G$-function $F(z)$. 

\medskip

The quality of restricted rational approximants as in Theorem \ref{th1} and Corollary \ref{coro} can be measured (when $t=0$) by a Diophantine exponent $v_b$ studied in \cite{AmB}. Given $\xi\in\R\moins\Q$, $v_b(\xi)$ is the infimum of the set of real numbers $\mu$ such that $|\xi - \frac{n}{b^m}| \geq b^{-m(1+\mu)}$ for any $n\in\Z$ and any sufficiently large $m$. With this notation, the special case $t=0$ of Corollary \ref{coro} reads $v_b(F(a/b))\leq \eps$. The metric properties of this Diophantine exponent are studied in \cite[Section 7]{AmB}: with respect to 
Lebesgue measure, almost all real numbers $\xi$ satisfy $v_b(\xi) = 0$ for any $b\geq 2$, and given $b\geq 2$ the set of $\xi$ such that $v_b(\xi) \ge \eps$ has Hausdorff dimension $\frac{1}{1+\eps}$. Therefore Theorem \ref{th1} and Corollary \ref{coro} are a 
step towards the conjecture that values of $G$-functions behave like generic real numbers with respect to rational approximation.

\medskip

Our results have interesting consequences on the   nature of the $\bz$-ary expansions of  values of $G$-functions; this is a class of numbers for which very few such results are known (see \cite{BK2, BK1}). 
Let $\bz$, $t$ be integers with $\bz\geq 2$ and $t\geq 1$, and let $\xi\in\R\setminus\Q$. We denote by $0.a_1a_2a_3\ldots$ the expansion in base $\bz$ of the fractional part of $\xi$. For any $n\geq 1$,  let $\notn$  denote the largest integer $\ell$ such that $(a_na_{n+1}\ldots a_{n+t-1})^\ell$ is a    prefix of the infinite word  $a_na_{n+1}a_{n+2}\ldots $.   In other words, it is the number of times the pattern $a_na_{n+1}\ldots a_{n+t-1}$ is repeated starting from $a_n$. Obviously $\notn \geq 1$, and $\notn$ is finite since $\xi$ is irrational. If $t=1$, $\notn$ is simply the number of consecutive equal digits in  the expansion of $\xi$, starting from $a_n$. For almost all real numbers $\xi$ with respect to Lebesgue measure, $\lim_{n\to\infty}\frac1{n}\,\notn =0$.

\begin{theo} \label{th2}
 Let $F$ be a $G$-function with rational Taylor coefficients and with $F(z)\not\in\Q(z)$, $\eps > 0$,  
 and $a\in \Z$, $a\neq 0$.
Let $\bz\geq 2$. Then for any $s\geq 1$ such that $\bz^s$ is 
sufficiently large in terms of $F$, $\eps$,   and $a$, we have  for any $t\geq 1$:  
$$\limsup_{n\to\infty}\frac1{n}\,\notnbis \leq 
\eps / t.$$
\end{theo} 
In the case of the dilogarithm, this result applies to $\textup{Li}_2(1/b^s)$ for   $a=1$,  any fixed $\eps \in (0,1)$, any $t\ge 1$ and any $b\ge 2$ provided $s\ge 10^{7}/\eps$.  
A similar  bound on  this upper limit, but with $1+\eps$ instead of $\eps$, follows from (and under the assumptions of) Theorem \ref{theo:zud}. Conjecturally, we have $\lim_n\frac1{n}\,\notn = 0$ whenever $\xi$ is a  transcendental value of a $G$-function, but it seems that the only such $\xi$ for which the upper bound  $\limsup_{n}\frac1{n}\,\notnun  < 1$ was known are values of the logarithm~\cite{riv1}.

\medskip

When $\xi$ is an irrational algebraic number, Ridout's theorem~\cite{ridout} yields $v_b(\xi) =0$  and $\lim_{n}\frac1{n}\,\notn = 0$ for any $b$ and any $t$.  It is not effective: for a general real algebraic  number $\xi$, given $b$, $t$ and $\eps>0$, no explicit value of $N(\xi,b,t,\eps)$  is known such that $\notn \leq \eps n$ for any $n \geq N(\xi,b,t,\eps)$. On the contrary, if $\xi = F(a/b^s)$ then  Theorem~\ref{th2} (which is effective) provides such an explicit value provided $b^s$ is large enough -- recall that algebraic functions which are holomorphic at 0 are $G$-functions. However if $\xi$ is fixed then Theorem \ref{th2} applies only if $\eps$ is not too small: we do not really get an effective version of Ridout's theorem for $\xi$. For other results concerning the $b$-ary expansions of algebraic numbers, we refer the reader to~\cite{ab, dbcp, riv3}.

We proved in~\cite{firi} that any real algebraic number is equal to $F(1)$ 
for some algebraic $G$-function $F$ with rational coefficients and 
radius of convergence arbitrarily large. Unfortunately, we do not have a control on 
the growth of 
the sequence of denominators of the coefficients of $F$, which is important in the computation of the constants in Theorem~\ref{th1}. Therefore, we cannot prove that any real algebraic number can be realized as a $G$-value 
$F(a/b)$ to which Theorem~\ref{th1} applies. 

\medskip

Finally, let us explain the basic reason behind our improvement on Zudilin's exponent. To estimate the difference 
$\vert F(\frac ab)-\frac pq\vert$
using the methods of this article, we need at some point to find a lower bound on  a certain difference 
$D=\vert \frac pq-\frac{u}{b^{k}v}\vert$ between two distinct 
rationals ($k\in\mathbb{N} = \{0,1,2,\ldots\}$, $p,q,u,v \in\mathbb{Z}$). 
When $q$ could be anything, 
the best we can say is that, trivially, $D \ge (b^{k}qv)^{-1}$; however, if we know 
in advance that $q=B\cdot b^{m}$ then we can improve the trivial bound to 
$D\ge (b^{\max(m,k)}Bv)^{-1}$  and we save a  factor $b^{\min(m,k)}$ in the process. The fraction $\frac{u}{b^{k}v}$ is obtained by constructing (inexplicit)  
Pad\'e approximants of type II to $F(z)$ and the other $G$-functions appearing in a differential system of order $1$ satisfied by~$F$. Inexplicit Pad\'e approximation is a classical tool in the Diophantine theory of $G$-functions. 

Our main new ingredient is the use of  non-diagonal Pad\'e type approximants, i.e. the polynomials are made to have different degrees, which creates the factor $b^{k}$ we need. 
This idea seems to have been introduced in~\cite{beukers1} in a particular case; we use it in its full generality. To illustrate its 
importance for Theorem~\ref{th2}, we remark that if one tries to compute  
an irrationality measure for $F(a/b)$ under the assumptions of 
Eq.~\eqref{eqirratf}  
with the method of the present paper, one gets an irrationality exponent not smaller than $N+1+\eps$, 
where $N$ is the least integer such that $1, F(z), \ldots, F^{(N)}(z)$ are linearly dependent over $\mathbb Q(z)$.

\medskip

The structure of this paper is as follows. Section \ref{subsec21} is devoted to general results on Pad\'e type approximation, and Section \ref{secpreuve} to the proof of Theorem \ref{th1}. At last, we deduce Theorem~\ref{th2} in Section \ref{secfac} and conclude with some remarks in  Section \ref{sec:conclusion}.

\section{Non-diagonal Pad\'e type approximation}  \label{subsec21}

We gather in this section known results and preparatory computations that will be used in 
Section \ref{secpreuve} to prove Theorem \ref{th1}.

\subsection{Setting and zero estimate}\label{ssec:zeroest}

Let $F_1(z), \ldots, F_N(z)$ be $G$-functions with rational coefficients. We let $F_0(z) = 1 $ and assume that $F_0, F_1,\ldots,F_N$ are linearly independent over $\Q(z)$. We assume also that  $Y(z)={}^t\big(1, F_1(z), \ldots, F_N(z)\big)$  is  a solution of a 
 differential system of order~$1$
\begin{equation}\label{eq:sysdiffn}
Y'(z)=A(z)Y(z)
\end{equation}
where  $A(z)\in M_{N+1,N+1}(\mathbb Q(z))$ is a matrix of which the rows and columns are numbered from $0$ to $N$.

\medskip

Let $\mathcal{D}(z)$ be a non-zero polynomial in $\mathbb Z[z]$    such that 
$\mathcal{D}(z)A(z)\in M_{N+1,N+1}(\mathbb Z[z])$.   Let $d\in\N$ be such that 
$$\deg \calD(z) \leq d \;\mbox{ and }\; \deg \calD(z) A_{i,j}(z) \leq d-1$$ for any coefficient $A_{i,j} (z)$ of $A(z)$. 
We observe that $\mathcal{D}(z)$ is not a constant polynomial because if $A(z)$ has polynomial entries, the system \eqref{eq:sysdiffn} cannot have a non-zero vector of solutions consisting of $G$-functions; therefore $d\geq 1$.

\medskip

For any integers $p,q,h$ such that $p\ge q\ge Nh\ge 0$, we can find $N$ polynomials $P_{1}(z), \ldots, P_{N}(z)\in \mathbb Q[z]$ of degree $\le p$ and $Q(z) \in \mathbb Q[z]$ of degree $\le q$,  such that the order at $z=0$ of 
$$
R_j(z):=Q(z)F_j(z)-P_j(z)
$$ 
is at least $p+h+1$ for all $j=1,\ldots, N$. In particular, $Q(z)$ is not identically zero. 
We say that $(Q;P_1, \ldots, P_N)$ is a Pad\'e type approximant of type  II $[q;p, \ldots, p;p+h+1]$ of $(F_1, \ldots, F_N)$. 
It is not unique in general. 

In what follows it is convenient to let $P_0(z) = Q(z)$ and $R_0(z) = 0$, even though they do  not play exactly the same role as the other $P_j$'s and $R_j$'s.

\medskip

Set ${\bf P}(z)={}^t\big(P_{0}(z), \ldots, P_{N}(z)\big)$ and ${\bf R}(z)={}^t\big(R_{0}(z), \ldots, R_{N}(z)\big)$. 
Following Chudnovsky~\cite{chud1, chud2}, 
for $k\ge 0$ we define ${\bf P}_{k}(z):={}^t\big(P_{0,k}(z), \ldots, P_{N,k}(z)\big) \in \mathbb Q[z]^{N+1}$ and   ${\bf R}_{k}(z):={}^t\big(R_{0,k}(z), \ldots, R_{N,k}(z)\big) \in \mathbb Q[[z]]^{N+1}$ by 
\begin{equation} \label{eqdefppp}
{\bf P}_{k}(z):=\frac1{k!}\mathcal{D}(z)^k\Big(\frac{d}{d z}-A(z)\Big)^{k}{\bf P}(z),
\end{equation}
$$
{\bf R}_{k}(z):=\frac1{k!}\mathcal{D}(z)^k\Big(\frac{d}{d z}-A(z)\Big)^{k}{\bf R}(z).
$$

Now recall that $F_0=1, F_1,\ldots,F_N$ are linearly independent, so that the matrix $A(z)$ is uniquely determined by these functions and the zero-th row of $A$ is identically zero. Therefore we obtain the formula
\begin{equation} \label{eqdefqqq}
Q_{k}(z) =\frac1{k!}\mathcal{D}(z)^kQ^{(k)}(z)
\end{equation}
where $Q_k(z) := P_{0,k}(z)$. An important property is that if $Q(z)\in \mathbb Z[z]$, then $Q_k(z)\in \mathbb Z[z]$ for any $k$   because $\frac1{k!}(z^j)^{(k)}=\binom{j}{k} z^{j-k}$. 
 Moreover we have for any $k\in\N$ and any $j$:
 $$\deg Q_k \leq q+(d-1)k \quad \mbox{ and } \quad\deg P_{j,k} \leq p+(d-1)k.$$

We shall make use of the following results.
Part  $(i)$ follows easily from  the bounds on the degrees of $Q_k$ and $P_{j,k}$ and the relation $R_{j,k} = Q_k F_j - P_{j,k}$: see \cite[\S 2]{chud2}. Part $(ii)$ is the difficult one: it is a refinement and correction  by Andr\'e~\cite[p. 115]{andre} of Chudnovsky's zero estimate~\cite{chud1, chud2}. The fact that $F_1$, \ldots, $F_N$ are $G$-functions is used only to make the constant in  $(ii)$ effective.

\begin{theo}[Chudnovsky, Andr\'e] \label{lem:1} Let $(Q;P_1, \ldots, P_N)$ be  a Pad\'e type approximant of type  II $[q;p, \ldots, p ;p+h+1]$ of $(F_1, \ldots, F_N)$; 
recall that $F_0 (z) = 1$, $F_1(z)$, \ldots, $F_N(z)$ are $\mathbb Q(z)$-linearly independent $G$-functions with rational coefficients.
 Then:

$(i)$ For any $k\ge 0$ such that $h\ge kd$, 
 $(Q_k;P_{1,k}, \ldots, P_{N,k})$ is a Pad\'e type approximant 
$$
[q+k(d-1);p+k(d-1),\ldots, p+k(d-1);p+h+1-k]
$$ 
of $(F_1, \ldots, F_N)$.

\medskip

$(ii)$ The determinant 
$$
\Delta_N(z):=\left \vert 
\begin{array}{ccc}
Q_{0}(z)&\cdots & Q_{N}(z)
\\
P_{1,0}(z)& \cdots & P_{1,N}(z)
\\
\vdots & \vdots & \vdots
\\
P_{N,0}(z)&\cdots & P_{N,N}(z)
\end{array}
\right \vert
$$
is not identically zero provided $h\geq\hz$, where $\hz$ is a positive constant, which depends only on $F_1$, \ldots, $F_N$ and can be computed explicitly.
\end{theo}
Let us deduce precisely this result from Andr\'e's theorem. Given distinct  integers $i,j \in \{1,\ldots,N\}$ we have 
\begin{align*}
\ord_0 (P_iF_j - P_j F_i) 
&= -\ord_0(Q)+\ord_0(P_i(QF_j-P_j)-P_j(QF_i-P_i))\\
&\geq   -\ord_0(Q)+ \min(\ord_0 (P_i), \ord_0(P_j)) + p+h+1 \\
&\geq   p+h+1
\end{align*}
because $\ord_0 (P_i) \geq \ord_0 (Q)$ for any $i$. Since we also have 
 $\ord_0(QF_j-P_j) \geq  p+h+1$, Andr\'e's zero estimate \cite[p. 115]{andre}  applies as soon as $h$ is greater than the constant he denotes 
by $c_0(\Lambda)$, that we call $h_0$ here.  Moreover $h_0$ is effective: see  \cite[Exercise 2, p. 126]{andre}.
For future reference, we notice that $(i)$ and $(ii)$ can be combined as soon as $h \geq \max(\hz, Nd)$; this will be the case below.

\subsection{Construction of the Pad\'e approximants}

Let us explain precisely now  the construction of the $P_j$'s and of $Q$. First of all, we set 
$$
F_j(z)=\sum_{n=0}^{\infty} f_{j,n} z^n.
$$ 
Since the $F_j$'s are $G$-functions, there exist a sequence of integers $d_n>0$ and a constant $D>0$ such that $d_n f_{j,n} \in \mathbb Z$ and $d_n\le D^{n+1}$, and also  a constant $C>0$ such that 
$\vert f_{j,n}\vert \le C^{n+1}$ for all $n\ge 0$ and all $j$. Let us write $P_j(z)=\sum_{n=0}^p u_{j,n} z^n$ for $1\leq j \leq N$ and  $Q(z)=\sum_{n=0}^q v_{n} z^n$. 
By definition of the $P_j$'s and of $Q$, we have the equations
\begin{equation}\label{eq:11}
\sum_{k=0}^q f_{j,n-k} v_k =0, \qquad n=p+1,\ldots, p+h,  \;j=1, \ldots, N
\end{equation}
and
\begin{equation}\label{eq:21}
\sum_{k=0}^{\min(n,q)} f_{j,n-k} v_k =u_{j,n}, \qquad n=0,\ldots, p, \;j=1, \ldots, N.
\end{equation}

Multiplying Eq.~\eqref{eq:11} by $d_{p+h}$, we obtain a system of $Nh$ equations in the $q+1$ unknowns $v_0, \ldots, v_q$, 
with integer coefficients 
$d_{p+h}f_{j,0}, \ldots, d_{p+h}f_{j,p+h}$ bounded in absolute value by $(CD)^{p+h+1}$. Since $q+1>Nh$, Siegel's lemma (see for instance  \cite[Lemma 11, Chapter~3, p. 102]{Shidlovski})  implies 
the existence of a non-zero vector 
of solutions $(v_0, \ldots, v_{q})\in \mathbb Z^{q+1}$ such that
\begin{equation}\label{eq:13}
\vert v_k\vert \le 1 + (q(CD)^{p+h+1})^{\frac{Nh}{q+1-Nh}}, \quad k=0, \ldots, q.
\end{equation}
From Eq.~\eqref{eq:21}, we see that $d_{p}P_j(z)\in \mathbb Z[z]$ for $j=1, \ldots, N$. 
Let $H(A)$ denote the maximum 
of the moduli of the coefficients of a polynomial $A(z)$ with real coefficients. Since $CD\ge 1$, Eq.~\eqref{eq:13} implies that 
\begin{equation}\label{eq:HQ}
H(Q) \le  2(q(CD)^{p+h+1})^{\frac{Nh}{q+1-Nh}}.
\end{equation}

\subsection{Properties of $Q_k$, $P_{j,k}$ and $R_{j,k}$}

In this section, we collect some informations we shall  use freely in the proof of Theorem~\ref{th1}.

From the estimate~\eqref{eq:HQ}, we can bound the coefficients of $Q_k(z)\in \mathbb Z[z]$ for any $k\ge 0$. Indeed, 
we shall prove that\footnote{The proof of this inequality in the published version of this paper contains a mistake, pointed out to us by
 Dimitri Le Meur and corrected in an erratum.}
\begin{equation}
H(Q_k) 
  \le 2^{2q+(d-1)k+1}H(\mathcal{D})^k(q(CD)^{p+h+1})^{\frac{Nh}{q+1-Nh}}. 
\label{eqhqk} 
\end{equation}
For any $A,B\in\C[X]$ we have
$$H(AB) \leq \min(1+\deg A, 1+\deg B) \, H(A) \, H(B)$$
so that 
$$H( \calD^k ) \leq (1+\deg\calD)^{k-1} H(\calD)^k \leq (d+1)^{k-1}  H(\calD)^k $$
and
$$H(Q_k) \leq c_k  H(\calD)^k H(Q)$$
for any $k\geq 0$, where $c_k=0$ if $k > q$ and, if $k\leq q$,
\begin{eqnarray*}
c_k &=& \min( 1+ k\deg\calD, 1+ \deg Q^{(k)})  (d+1)^{k-1}  2^q \\
&\leq & (q-k) (d+1)^{k-1}  2^q \\
&\leq & 2^{q-k}  (d+1)^{k}  2^q\\
&\leq & 2^{2q}  \Big(\frac{d+1}2\Big)^{k}  \\
&\leq &2^{2q+(d-1)k}
\end{eqnarray*}
where the last inequality comes from the fact that $\frac{d+1}2 \leq 2^{d-1}$ for any positive integer $d$. Taking Eq. \eqref{eq:HQ} into account, this concludes the proof of Eq. \eqref{eqhqk}.

Let us  now bound $R_{j,k}(z)$ for $1\leq j \leq N$.  
Letting $Q_k(z) = \sum_{n=0}^{q+(d-1)k}v_n^{(k)}z^n$ and recalling that $R_{j,k} = Q_k F_j - P_{j,k}$, Lemma \ref{lem:1} $(i)$ yields 
$$
R_{j,k}(z)=\sum_{n=p+h+1-k}^{\infty} \Big(\sum_{\ell=0}^{\min(n, q+k(d-1))} f_{j,n-\ell} v_{\ell}^{(k)}\Big)z^n
$$
from which we deduce that, for $|z| < 1/C$:  
\begin{align}
\vert R_{j,k}(z)\vert &\le H(Q_k)(q+k(d-1)+1) \max(1,C)^{q+k(d-1)} \sum_{n=p+h+1-k}^{\infty} C^n  \vert z\vert ^n\nonumber
\\
&\le \frac{ H(Q_k)(q+k(d-1)+1)}{1-C\vert z\vert}  \max(1,C)^{q+k(d-1)} (C\vert z\vert)^{p+h+1-k}. 
\label{eqrjk} 
\end{align}

Finally, for $j=1, \ldots, N$, letting $P_{j,k}(z) = \sum_{n=0}^{p+(d-1)k}u_{j,n}^{(k)}z^n $  
we have 
\begin{equation*}
\sum_{m=0}^{\min(n,q+(d-1)k)} f_{j,n-m} v_m^{(k)} =u_{j,n}^{(k)}, \qquad n=0,\ldots, p+(d-1)k.
\end{equation*}
It follows that $d_{p+(d-1)k}P_{j,k}(z)\in \mathbb Z[z]$ for all $k\ge 0$ and $1\le j \le N$.

\section{Proof of Theorem~\ref{th1}} \label{secpreuve}

We split the proof into two parts: in Section \ref{subsec22}  we prove a general (and technical) result, namely Eq. \eqref{eqccltech}, from which Theorem \ref{th1} will be deduced in Section~\ref{subsec23}.

\subsection{Main part of the proof} \label{subsec22}

 We keep the notation and assumptions of Section~\ref{subsec21} concerning $F_1$, \ldots, $F_N$, $A$, $\calD$, $d$, $C$, $D$.  Without loss of generality, we may assume that $C\geq 1$.

We fix   $t,x,y \in\R$ and $ a,b,B,m,n,h\in\Z$ such that $b,m\geq 1$ and 
$$
1 \leq \left|\frac ab\right| <  \min\left(\frac{1}{2C}, \frac1{H(\calD)}\right),
\hspace{0.5cm}
1\leq B \leq b^t, 
\hspace{0.5cm}
 0 < y < \frac1d , 
 \hspace{0.5cm}
  x > N+y.$$
We also assume that $h$ is sufficiently large; in precise terms, we assume that 
$$ 
h\geq \max\left(\hz, 8 N^2d^3, 
\frac{m}{x-N-y}\right) 
$$
where $\hz$ is the constant in Lemma \ref{lem:1},
and we shall also assume below (just before Eq. \eqref{eqhyp}) that $h$ is greater than some other positive constant that could be effectively computed in terms of $F_1$, \ldots, $F_N$.

We let $\beta = b^{t/h}$ and make one more assumption on these parameters, namely Eq. \eqref{eqhyp} below. At last, we fix an integer $j\in\{1,\ldots,N\}$. Then we shall deduce a lower bound on 
$\left| F_{j}(\frac ab) -\frac{n}{B\cdot b^m} \right| $, namely Eq. \eqref{eqccltech}.

\bigskip

Changing $z$ to $-z$ in all functions $F_1$, \ldots, $F_N$, we may assume that $a >0$. 

Let $z_0 \neq 0$ be a rational root of $\calD$; let us write $z_0 = r_0/r_1$ with coprime integers $r_0$, $r_1$. Then $r_1$ divides the leading coefficient of $\calD$, so that $|r_1 | \leq H(\calD)$ and $|z_0| \geq \frac1{|r_1|} \geq \frac1{H(\calD)}$. Therefore $a/b$ is not a root of $\calD$.

To apply the constructions of Section~\ref{subsec21} we let 
$$p = \pe{xh} \hspace{0.5cm}
 \mbox{ and }
 \hspace{0.5cm}
 q = \pe{(N+y)h},$$
  so that 
$$p \geq q+m.$$

Let us choose $k$ now. The determinant $\Delta_N(z)$ of Lemma \ref{lem:1} $(ii)$ has degree at most $q+Np + (d-1)N(N+1)/2$. We use the vanishing properties of Lemma \ref{lem:1} $(i)$ (since $h\geq Nd$) by susbtracting $F_i(z)$ times the zero-th row from the $i$-th row, for any $1\leq i \leq N$. We obtain that  $\Delta_N(z)$ vanishes at 0 with multiplicity at least $N (p+h+1)-N(N+1)/2$. Therefore we have
$$
\Delta_N(z)=z^{ N (p+h+1)-\frac{N(N+1)}2} \widetilde{\Delta}_N(z)
$$
where $\widetilde{\Delta}_N(z)$ has degree $\le \ell_0$, with $\ell_0 = q  - N (h+1) +dN(N+1)/2$,  and is not identically zero.
Since $a/b$ is different from 0, its multiplicity as a root of $ \Delta_N(z)$ is at most $\ell_0$. Following the proof of \cite[Theorem 4.1]{chud2},   we deduce that the matrix
$$
\left(
\begin{array}{ccc}
Q_{0}(a/b)&\cdots & Q_{N+\ell_0}(a/b)
\\
P_{1,0}(a/b)& \cdots & P_{1,N+\ell_0}(a/b)
\\
\vdots & \vdots & \vdots
\\
P_{N,0}(a/b)&\cdots & P_{N,N+\ell_0}(a/b)
\end{array}\right)
$$
has rank $N+1$. Therefore we have $n Q_k(a/b)-Bb^m P_{j,k}(a/b)\neq 0$ for some integer $k$, with
$$k \leq \ell_0+N =   q - N h +dN(N+1)/2 ; $$
recall that $j$ is fixed in this proof.

\bigskip

 By construction of the polynomials $P_{j,k}(z)$ and $Q_k(z)$, there exist 
two integers $U_{j,k},V_k$  such that $P_{j,k}(a/b)=U_{j,k}/(d_{p+(d-1)k} b^{p+(d-1)k})$ and $Q_k(a/b)=V_k/b^{q+(d-1)k}$.
We deduce that  
$$\xi = d_{p+(d-1)k} b^{p+(d-1)k} \big( n Q_k(a/b) - Bb^m P_{j,k}(a/b) \big)$$
is a non-zero integer (since $p\geq q$). Moreover we have assumed that $p\geq q+m$ so that $\xi$ is a multiple of $b^m$, and thus 
$$|\xi| \geq b^m.$$

On the other hand we have
$$\xi = d_{p+(d-1)k} b^{p+(d-1)k} \Big(  Q_k(a/b) \big(n- Bb^m F_{j}(a/b) \big) - Bb^m \big(  P_{j,k}(a/b) - Q_k(a/b) F_{j}(a/b) \big) \Big)$$
so that
$$|\xi| \leq  d_{p+(d-1)k} b^{p+(d-1)k} \Big( \big|Q_k(a/b)\big| \cdot \big|n- Bb^m F_{j}(a/b) \big| + Bb^m  \big| R_{j,k}(a/b)\big|\Big).$$
Comparing this upper bound and the lower bound $|\xi| \geq b^m $ we obtain
\begin{equation} \label{eq16}
\big|Q_k(a/b)\big| \cdot \big|n- Bb^m F_{j}(a/b) \big|  \geq  d_{p+(d-1)k}^{-1} b^{-p-(d-1)k+m} -  Bb^m  \big| R_{j,k}(a/b)\big|.
\end{equation}
We shall prove below that under a suitable assumption (namely Eq. \eqref{eqhyp}) we have
\begin{equation} \label{eqRnegli}
 \big| R_{j,k}(a/b)\big| < \frac12 d_{p+(d-1)k}^{-1} b^{-p-(d-1)k }B^{-1}
\end{equation}
so that  the right hand-side of Eq. \eqref{eq16} is positive, and $Q_k(a/b) \neq 0$. Moreover  Eq. \eqref{eq16} yields
\begin{equation} \label{eqminotemp}
\left| F_{j}\Big(\frac ab\Big) -\frac{n}{B\cdot b^m} \right| \geq \frac{d_{p+(d-1)k}^{-1} b^{-p-(d-1)k }  }{2B |Q_k(a/b)|}.
\end{equation}

Now, recall that 
$$
p =\pe{xh} , \mbox{ } q = \pe{( N+y)h}, \mbox{ and  }  k \leq yh +  \frac{dN(N+1)}2.
$$
Let us denote by  $\gdo(1)$ any positive quantity that can be bounded (explicitly) in terms of $F_1$, \ldots, $F_N$; such a bound may involve, among others, $d$, $N$, $\calD$, $C$ or $D$. We recall that $C,D\geq 1$ and notice that $\frac{Nh}{q+1-Nh} \leq \frac Ny$. Then Eq. \eqref{eqhqk}  yields 
$$H(Q_k) \leq \big( 2^{2(N+y)+(d-1)y } H(\calD)^y (CD)^{(x+1)N/y}\big) ^h \cdot (CD(N+y)h)^{N/y} \cdot  \gdo(1) $$
so that Eq. \eqref{eqrjk} provides, since $Ca/b<1/2$: 
\begin{multline*}
\big| R_{j,k}(a/b)\big| \leq 
\\
\big( 2^{2(N+y)+(d-1)y } H(\calD)^y (CD)^{(x+1)N/y}  C^{N+dy} (C a/b)^{x+1-y} \big) ^h \cdot (CD (N+dy)h)^{1+N/y}\cdot  \gdo(1).
\end{multline*}
Now let us assume, for simplicity, that $y \geq \frac1{8d}$. Then we have $(d-1)k \leq dyh$ since $hy\geq  N^2d^2$, so that
$$
d_{p+(d-1)k}  b^{ p+(d-1)k }B \leq \big(\beta (bD)^{x+dy} \big)^{h }.
$$
Therefore \eqref{eqRnegli} holds if $h$ is larger than some effectively computable constant (depending only on  $F_1$, \ldots, $F_N$)  and  if
\begin{equation} \label{eqhyp}
y \geq \frac1{8d} \;\mbox{ and }\;
2^{2(N+y)+(d-1)y } H(\calD)^y (CD)^{(x+1)N/y}  C ^{N+dy} (C  a/b )^{x+1-y}  (bD)^{x+ d  y} \beta < \frac12.
\end{equation}
Moreover, since $ a/b <1$, we have
\begin{eqnarray*}
\big |Q_k(a/b)\big| &\leq& (q+(d-1)k+1) H(Q_k) \\
 &\leq& \Big( 2^{2(N+y)+(d-1)y } H(\calD)^y (CD)^{(x+1)N/y}\Big) ^h \cdot  (CD (N+dy)h)^{1+N/y}\cdot  \gdo(1)
\end{eqnarray*}
so that Eq. \eqref{eqminotemp} yields finally (since $1/y \leq 8d$):
\begin{multline} \label{eqccltech}
 \left| F_{j}\Big(\frac ab\Big) -\frac{n}{B\cdot b^m}\right|  \ge
\\
\Big( \beta^{-1} (bD)^{-x- d y} 2^{-2(N+y)-(d-1)y } H(\calD)^{-y} (CD)^{-(x+1)N/y} \Big)^{h} \cdot h^{-9Nd} \cdot \gdo(1)^{-1}.
\end{multline}
This is a very general lower bound and in the next section we will proceed to a suitable choice of the parameters.

\subsection{Choice of the parameters and conclusion} \label{subsec23}

In this section we prove Theorem \ref{th1} by applying the proof of Section~\ref{subsec22} to suitable parameters.

  Let $F$ be a   $G$-function with $F(z)\not\in\Q(z)$. Let $F_0(z) = 1$, and denote by $N$ the least positive integer for which there exist $  a_0(z), \ldots, a_N(z)\in\Q(z)$ such that 
  $$F^{(N)}(z) = a_N(z) F^{(N-1)}(z) + \ldots + a_2(z) F'(z) + a_1(z) F(z)  + a_0(z).$$
We have $N \geq 1$, and $N=1$ may hold (it does for instance with $F(z) = \log(1-z)$). By construction and since $F(z)\not\in\Q(z)$, the $G$-functions $F_0 = 1$, $F_1 = F$, $F_2 = F'$, \dots, $F_N = F^{(N-1)}$ are linearly independent over $\Q(z)$. We are in position to apply the results of Sections~\ref{subsec21} and \ref{subsec22}; as in Section~\ref{subsec22} we may assume $a>0$ and $C\geq 1$.

\medskip

We let 
$$
\ch = 4 H(\calD) (CD)^{8Nd+1}  C   \mbox{ and } c_2 = 3(N+2). 
$$
We take $\cs = \max(\hz, \hun, \hd, 8N^2d^3,  4t)$  where $\hz$ is the constant implied in Lemma \ref{lem:1} $(ii)$,   $\hun$ is the effectively computable constant defined just before Eq. \eqref{eqhyp}, and $\hd  $ is another effectively computable constant to be defined below; these three constants depend only on $F$. Then we assume
\begin{equation}\label{eqminom}
m \geq \frac{\cs}3 \frac{\log(b)}{\log(\ch a )} ; 
\end{equation}
this is a consequence of our assumption   $m\geq c_3 \frac{\log(b)}{\log(a+1)}$ provided we choose $c_3 = \frac13 c_5$.

We choose
$$y=\frac1{4(d+1)}, \hspace{1cm} x = \frac13  \frac{\log(b)}{\log(\ch a )} \hspace{1cm} \mbox{ and } \hspace{1cm} 
h =\left\lfloor{\frac{m}{x-N-1}}\right\rfloor.
$$
Then \eqref{eqminob} implies $x\geq N+2$, and  \eqref{eqminom} yields $h\geq \cs $. 

Let us check that \eqref{eqhyp} holds. We notice that 
 $\beta   \leq   b^{1/4}$ 
since $h \geq \cs \geq 4t$. Since $y = \frac1{4(d+1)} $ and $x\geq N+1$, the left hand side of  \eqref{eqhyp} is less than
$$ 2^{2N+1}H(\calD)  (CD)^{8Nd(x+1)}   C^{x+N+2 }D^{x+1}a^{x+1}b^{-1/2}  
 \leq  (\ch a)^{x+1}b^{-1/2} \leq 1/2; 
$$
indeed  we have used the definition of $x$ and the lower bound $\frac{\log(b)}{\log(\ch a )}  \geq 9 $ which follows from  $x\geq N+2\geq 3$.
Therefore  all the assumptions made in  Section~\ref{subsec22} hold.

We set 
\begin{equation*}
\cn =D^{1+\frac{d}{(N+2)(d+1)}}2^{\frac{8N+1}{4N+8}}H(\mathcal{D})^{\frac{1}{4(N+2)(d+1)}}(CD)^{\frac{4N(N+3)(d+1)}{N+2}}.
\end{equation*} 
Using \eqref{eqccltech} and the various conditions on $x$ and $y$, we readily obtain
\begin{equation} \label{eq:rajout}
\left\vert F \Big(\frac ab\Big) - \frac{n}{B\cdot b^m}\right\vert \ge 
\beta^{-h} \Big( b^{x+dy} \cn^x\Big)^{-h} \cdot h^{-9Nd}\cdot \gdo(1)^{-1} \geq 
  B^{-1} (b \cn)^{-\frac{x+1}{x-N-1}\cdot m}
\end{equation}
provided $h\geq \hd$, where $\hd$ is effective and depends only on $F$. Now we have
$$
\frac{x+1}{x-N-1}=1+\frac{N+2}{x}\cdot\frac{1}{1-\frac{N+1}x}
\le 1+3(N+2)^2\frac{\log(c_1 a)}{\log(b)}
$$
because $x=\frac13\frac{\log(b)}{\log(c_1 a)}\ge N+2$. It is trivial matter to check that for any $u\ge1$ and any $v\ge e$, we have 
$\log(uv)\le 2\log(u+1)\log(v)$. Since $c_1\ge 4>e$ (because we always have $H(\mathcal{D})\ge 1$, $C\ge 1$ 
and $D\ge 1$), we can apply this with $u=a$, $v=c_1$ and we get 
$$\frac{x+1}{x-N-1} \leq 1  + \cdix \frac{\log(a+1)}{\log b}$$
with $c_{7}=6(N+2)^2\log(c_1)$. Hence, we deduce from \eqref{eq:rajout} that
$$
\left\vert F \Big(\frac ab\Big) - \frac{n}{B\cdot b^m}\right\vert \ge \frac1{
 B\cdot  b^{m} (a+1)^{c_{8} m} \cn^{m} } \ge  \frac1{
 B\cdot  b^{m} (a+1)^{c_4 m}}
 $$
where 
$$
c_{8}=\frac{\log(2c_6)}{\log(2)}c_{7} \quad \textup{and } \quad  
c_4=c_{8}+\frac{\log(c_6)}{\log(2)}.
$$
This completes the proof of Theorem \ref{th1}.

\medskip

\begin{Remark} 
Let us compute the constants $c_1, c_2$ and $c_4$ in the case  of the $G$-function $\textup{Li}_2$. The vector $Y(z)={}^t\big(1,\textup{Li}_1(z), \textup{Li}_2(z)\big)$ is solution of the differential system 
$$
Y'(z)= 
\left(
\begin{array}{ccc}
0 & 0 & 0\\ \frac{1}{1-z} & 0 & 0 \\ 0 & \frac{1}{z}& 0 \end{array}
\right)Y(z).
$$
Hence $\mathcal{D}(z)=z(1-z)$, $H(\mathcal{D})=1$, $d=2$, $N=2$, $C=1$ and $D=e^2$. 
With the constants   defined in the proof just above, we obtain $c_1=4e^{66}$, $c_2=12$ and 
$$
c_4=\frac{1201779}{48}+\frac{1185019}{3\log(2)}+396\log(2)< 10^{5.78}.
$$
The constant $c_3$ could be computed as well, but we did not try to do so because it is 
not important for the application to Theorem~\ref{th2}.

It follows that Theorem~\ref{th1} can be applied with $F(z)=\textup{Li}_2(z)$ when $b\ge e^{809} \vert a\vert^{12}>0$. 
Furthermore, Theorem~\ref{th2} can be applied for any integer $s$ such that $b^s\ge e^{809} \vert a\vert^{12}>0$ and
$b^s\ge (\vert a\vert+1)^{2c_4/ \eps }$. In particular, if   $a=1$, $b\ge2$ and $0<\eps<1$, Theorem~\ref{th2} applies to $\textup{Li}_2(1/b^s)$ for any integer $s \geq 10^{6.08}/\eps.$ We have not tried to optimize our general 
constants which in this case could be decreased. 
\end{Remark}

\section{Proof of Theorem~\ref{th2}} \label{secfac}

In this section we deduce Theorem~\ref{th2} from Corollary \ref{coro} stated in the introduction.

\medskip

Let $\xi = F(a/\bz^s)$,  $q_n  = \bz^{n-1}(\bz^t-1)$,  and 
$$
p_n = (\bz^t-1)\pe{\bz^{n-1}\xi} + a_n \bz^{t-1} + a_{n-1}\bz^{t-2} + \ldots + a_{n+t-1}.
$$ 
Then   the $\bz$-ary expansion of 
$$\frac{p_n}{q_n} = \frac{\pe{\bz^{n-1}\xi}}{\bz^{n-1}} + \frac{ a_n \bz^{t-1} + a_{n-1}\bz^{t-2} + \ldots + a_{n+t-1} }{\bz^{n-1} (\bz^t-1)}$$
has the same $n+t \notn-1$ first digits as the $\bz$-ary expansion  of $\xi$. Therefore we have
$$
\left|\,\xi - \frac{p_n}{q_n}
\right| \leq \frac{\bz-1}{\bz^{n+t\notn}}.
$$
Now Corollary \ref{coro} with 
$b^s$ for $b$, $B=b^t-1$ and $m = \pe{\frac{n-1}s}$ yields
$$
\left|\,\xi - \frac{p_n}{q_n}
\right| \ge \frac{1}{b^{\pe{\frac{n-1}s}s(1+\eps)}}.
$$ 
The comparison of both inequalities enables us to conclude the proof.

\section{Concluding remarks} \label{sec:conclusion}

In Section~\ref{subsec21}, we assumed that the degrees of the polynomials satisfy 
$p\ge q$ and in fact $p\ge q+m$, which was crucial to prove Theorem~\ref{th1}. The case $q\ge p$ also provides some informations, but not in the exact situation of Theorem~\ref{lem:1}. Indeed, with the notation of Section~\ref{ssec:zeroest} the polynomials $P_{j,k}(z)$ with $1\leq j \leq N$ depend on $P_1(z)$, \ldots, $P_N(z)$ and also on $P_0(z) = Q(z)$ (see Eq. \eqref{eqdefppp}). In order to be able to bound the degree of $P_{j,k}(z)$ in terms of $p$ only (independently of $q$), we need to deduce from \eqref{eqdefppp} a relation analogous to \eqref{eqdefqqq}, namely an expression for the polynomials $P_{j,k}(z)$  in terms of $P_1(z)$, \ldots, $P_N(z)$ only. This follows easily under the additional assumption that the zero-th row of $A(z)$ is identically zero, i.e. that ${}^t(F_1 , \ldots, F_N)$ is a solution of a  homogeneous linear differential system.
 Following the same method as in the 
case $p\ge q$, this enables us to prove the following result.
\begin{theo} \label{th11}
 Let $F$ be a  $G$-function with rational Taylor coefficients and $t\geq0$. Let us assume that $F(z)$ is solution of a homogeneous linear differential equation of order $N$ with coefficients in $\mathbb Q(z)$ and that $1, F(z), F'(z), \ldots, F^{(N-1)}(z)$ are linearly independent over $\mathbb Q(z)$. 
Then there exist some positive effectively computable constants $\widetilde{c}_1$, $\widetilde{c}_2$, $\widetilde{c}_3$, $\widetilde{c}_4$, depending only on $F$ (and $t$ as well for $\widetilde{c}_3$), such that the following property holds. 
Let $a\neq 0$ and $b,B\geq 1$ be integers such that 
\begin{equation}
b > (\widetilde{c}_1 |a|)^{\widetilde{c}_2} \mbox{ and } B \leq b^t. \label{eqminob2}
\end{equation}
Then $F(a/b)\notin \mathbb Q$ and for any $n\in\Z$ and any $m\geq \widetilde{c}_3 \frac{\log(b)}{\log(|a|+1)}$, we have  
\begin{equation}\label{eq:mesres2}
\left\vert \frac1{F (\frac ab)} - \frac{n}{B\cdot b^m}\right\vert \ge \frac{1}{B\cdot b^m \cdot (|a|+1)^{\widetilde{c}_4m} }.
\end{equation}
\end{theo}
Analogues of Corollary~\ref{coro} and Theorem~\ref{th2} for $1/F(a/b)$ hold as well. These results can be applied directly to the functions $\log(1-z)+\sqrt{1-z}$  and  $\sqrt{1-z}\log(1-z)$   for instance,  but not to $\log(1-z)$. Actually the proof of Theorem \ref{th11} (and of all other results in this paper) can be generalized to number fields, at least to multiply $B$ with a fixed non-zero algebraic number (and all implied constants would depend on this number), by replacing the algebraic number $\xi$ defined in Section \ref{subsec22} with its norm over the rationals. Applying Theorem \ref{th11} to  $\sqrt{1-z}\log(1-z)$  with $B$ multiplied by $\sqrt{1-a/b}$ and canceling out this factor shows that Theorem \ref{th1},  Corollary~\ref{coro} and Theorem~\ref{th2} hold with $1/\log(1-a/b)$ instead of $F(a/b)$.

\medskip

A natural problem is to obtain an analogue of Theorem~\ref{th1} when the  $F_j$'s are $E$-functions and not $G$-functions. With the same notations as in Section~\ref{subsec21}, the polynomials $Q_k$ would still have integer coefficients, but the denominators of the coefficients of the 
polynomials $P_{j,k}$ would no longer be bounded by $d_{p+(d-1)k}$ but by $(p+(d-1)k)!d_{p+(d-1)k}$. As the reader may check, this cancels the benefits of having non-diagonal Pad\'e type approximants if we follow the same method of proof as in Section~\ref{secpreuve}. We don't know if this problem can be fixed to prove analogues of Theorems~\ref{th1} and \ref{th2} for $E$-functions. Very few results are known on   $b$-ary expansions of values of $E$-functions (see  \cite{Boris}, \cite{BK2}, \cite{BK1}). From a conjectural point of view, the situation is not clear either: values of $E$-functions {\em do  not} behave like generic numbers with respect to rational approximation, as the continued fraction expansion of $e$ shows.

\def\refname{Bibliography}

\bigskip

\noindent S. Fischler, 
Laboratoire de Math\'ematiques d'Orsay, Univ. Paris-Sud, CNRS, Universit\'e Paris-Saclay, 91405 Orsay, France.

\medskip

\noindent 
T. Rivoal,  Institut Fourier,  CNRS et Universit\'e Grenoble Alpes, 
 CS 40700, 38058 Grenoble cedex 9, France

\medskip

\noindent Keywords: $G$-functions, rational approximations, irrationality measure, Pad\'e approximation, integer base expansion.

\medskip
 
\noindent MSC 2000: 11J82 (Primary); 11A63, 11J25, 11J91 (Secondary).

\pagebreak

\section*{Addendum}

{\em This addendum is not  included in the published version of this paper.}

\medskip

When $F(\frac ab)$ is an algebraic irrational, Corollary \ref{coro}  looks like Ridout's Theorem for  algebraic irrational numbers, but this is not really the same. First, if $F(\frac ab)$ is an algebraic irrational and $b$ is fixed, then it applies only if $\eps$ is not too small with respect to $b$, and thus we do not get an effective version of Ridout's theorem ``in base $b$'' for this number. Second, we don't know if any 
algebraic irrational number can be represented as a value $F(\frac ab)$ to which these results apply.

\medskip

In this addendum, we deduce from Theorem~\ref{th1} the following result, which partially solves these problems. 
\begin{theo}\label{thm:2} Let $d$ be a positive rational number   such that $\sqrt{d}\notin\mathbb Q$. There exist some constants $\eta_d>0, \kappa_d>0$ and $N_d$ such that
for any convergent $\frac{\alpha}{\beta}$ of the continued fraction expansion of  $\sqrt{d}$ with $\alpha, \beta\ge N_d$, we have 
$$
\left\vert \sqrt{d}-\frac{n}{\alpha^m}\right\vert \ge \frac{1}{(\eta_d\alpha)^{m}} \quad \textup{and} \quad \left\vert \sqrt{d}-\frac{n}{\beta^m}\right\vert \ge \frac{1}{(\kappa_d \beta)^{m}}
$$
for any integer $n\in \mathbb Z$ and any  $m$ large enough with respect to $d, \alpha, \beta$.
\end{theo}

In particular, for any $\eps> 0$ we have
$$
\left\vert \sqrt{d}-\frac{n}{\alpha^m}\right\vert \ge \frac{1}{ \alpha ^{m(1+\eps)}} \quad \textup{and} \quad \left\vert \sqrt{d}-\frac{n}{\beta^m}\right\vert \ge \frac{1}{  \beta ^{m(1+\eps)}}
$$
provided $\alpha$ and $ \beta$ are large enough (in terms of $d$ and $\eps$).

\begin{proof}  Let $\alpha, \beta$ be any positive integers such that $\vert \alpha^2-d\beta^2\vert \le c(d)$ for some given constant $c(d)$. 
Note that if $\alpha/\beta$ is a convergent to $\sqrt{d}$, then 
$$
\vert \alpha^2-d\beta^2\vert \le \frac{\alpha+\sqrt{d} \beta}{\beta} \le 2\sqrt d + 1 
$$
so that $c(d) = 2\sqrt d + 1 $ is an admissible value  for all convergents.

\medskip

Let $f(x)=\sqrt{1-x}$. Then 
$
f(\frac{\alpha^2-d\beta^2}{\alpha^2}) = \frac{\beta}{\alpha}\sqrt{d}.$ Let $d=\frac uv$ with positive integers $u$ and $v$. 
We can apply Theorem~\ref{th1} to $F =f$,   $a=v\alpha^2-u\beta^2$ and $b=v\alpha^2$,  provided 
that $\alpha^2>c_1^{c_2}\vert \alpha^2-d\beta^2\vert^{c_2}$ where $c_1, c_2$ depend only on $d$. This inequality holds a 
fortiori if we assume that  $\alpha\ge (c_1c(d))^{c_2/2}=:N_d$, which we now do. Then 
$$
\left\vert \frac{\beta}{\alpha}\sqrt{d}-\frac{n}{B \cdot(v\alpha^2)^{m}}\right\vert 
= \left\vert f\left(\frac{\alpha^2-d\beta^2}{\alpha^2}\right)  -\frac{n}{B\cdot (v\alpha^2)^{m}}\right\vert
\ge \frac{1}{B\cdot (1+v \vert \alpha^2-d \beta^2 \vert)^{c_4 m} \cdot (v\alpha^2) ^{m}}
$$ 
for any $1\le B\le v \alpha^{2t}$, any $n\in \mathbb Z$ and any $m\ge c_3 \frac{\log(v\alpha^2)}{\log(1+v \vert \alpha^2-d \beta^2 \vert)}$.

Thus
$$
\left\vert \sqrt{d}-\frac{\alpha n}{\beta B \cdot v^m \alpha^{2m}}\right\vert \ge \frac{\alpha}{\beta\cdot B\cdot (1+v c(d) )^{c_4 m} \cdot (v\alpha^2)^{m}}.
$$
Note that $c_3$ depends on $f$ and $t$. We now choose $t=2$, so that $c_3$ becomes absolute. With $B=\alpha$ and $n=\beta v^m  n'$ (for any $n'\in \mathbb  Z$), we get 
$$
\left\vert \sqrt{d}-\frac{n'}{\alpha^{2m}}\right\vert \ge \frac{1}{\beta\cdot (1+v c(d) )^{c_4 m} \cdot v^m \cdot \alpha^{2m}}.
$$
On the other hand, with $B=\alpha^2$ and $n=\beta v^m   n'$ (for any $n'\in \mathbb  Z$), we obtain
$$
\left\vert \sqrt{d}-\frac{n'}{\alpha^{2m+1}}\right\vert \ge \frac{1}{\beta\cdot (1+v c(d) )^{c_4 m} \cdot v^m\cdot \alpha^{2m+1}}.
$$
Moreover,  assuming $m\ge C(d,\alpha, \beta)$ we have 
$$
\beta\cdot (1+v c(d))^{c_4 m} v^m \le \delta^m
$$
for some constant $\delta$ that depends only on $d$. Therefore combining the previous inequalities yields
 $$
\left\vert \sqrt{d}-\frac{n'}{\alpha^m}\right\vert \ge \frac{1}{(\eta_d\alpha)^{m}}
$$
for any $n'\in\Z$ and any $m\ge C(d,\alpha, \beta)$, where 
 $\eta_d>0$ depends only on $d$.

\medskip

We now prove the other inequality
$$
\left\vert \sqrt{d}-\frac{n}{\beta^m}\right\vert \ge \frac{1}{(\kappa_d \beta)^{m}}.
$$
Any   convergent of $1/\sqrt{d}$ (except maybe the first ones)  is of the form $\beta/\alpha$ where $\alpha/\beta$ is a convergent of $\sqrt{d}$. Therefore we may apply the above result with $1/d$ and $\beta/\alpha$: we obtain 
 $$
\left\vert \frac1{\sqrt{d}}-\frac{n'}{\beta^m}\right\vert \ge \frac{1}{(\eta_d\beta)^{m}}.
$$
Since the map $x\mapsto 1/x$ is Lipschitz around $\sqrt d$, we deduce the lower bound of Theorem \ref{thm:2} by choosing an appropriate constant $\kappa_d$.
\end{proof}


\begin{thebibliography}{1}\label{sec:biblio}
\addcontentsline{toc}{section}{Bibliographie}

\bibitem{Boris} B. Adamczewski, {\em On the expansion of some exponential periods in an integer base}, Math. Ann. {\bf 346} (2010), 107--116.
\bibitem{ab} B. Adamczewski, Y. Bugeaud, {\em 
On the complexity of algebraic numbers I: Expansions in integer bases}, 
Ann. Math. (2) {\bf 165}.2 (2007), 547--565.
\bibitem{AmB} M. Amou, Y. Bugeaud, {\em Exponents of Diophantine approximation and expansions in integer bases}, 
 J. London Math. Soc. {\bf 81} (2010),  297--316.
 \bibitem{andre} Y. Andr\'e, {\em G-functions and Geometry}, Aspects of Mathematics, {\bf E13}. Friedr. Vieweg \&
Sohn, Braunschweig, 1989.
\bibitem{andrecrelle} Y. Andr\'e, {\em $G$-fonctions et transcendance}, J. reine angew. Math. {\bf 476} (1996), 95--125.
\bibitem{dbcp} D. H. Bailey, J. M. Borwein, R. Crandall and C. Pomerance, {\em On the binary expansions of algebraic numbers}, J. Th\'eor. Nombres Bordeaux {\bf 16}.3 (2004), 487--518.
\bibitem{BauerBennett} M. Bauer, M. Bennett, {\em Application of the hypergeometric method to the generalized Ramanujan-Nagell equation}, Ramanujan J. {\bf 6} (2002), 209--270.
\bibitem{BugeaudBennett} M. Bennett,  Y. Bugeaud, {\em Effective results for restricted rational approximation to quadratic irrationals}, Acta Arith. {\bf 155} (2012), 259--269. 
 \bibitem{BB} D. Bertrand,  F. Beukers,  {\em Equations diff\'erentielles lin\'eaires et majorations de multiplicit\'es}, Ann. Scient. Ec. Norm. Sup. {\bf 18} (1985),  181--192.
 \bibitem{beukers1} F. Beukers, {\em 
On the generalized Ramanujan-Nagell equation  I}, Acta Arith. {\bf 38}.4 (1980), 389--410.
 \bibitem{beukers} F. Beukers,  {\em A refined version of the Siegel-Shidlovskii theorem}, 
Ann. of Math. {\bf 163}.1 (2006),   369--379.
\bibitem{bombieri} E. Bombieri, {\em On G-functions}, Recent Progress in Analytic Number Theory {\bf 2}, Acad.
Press (1981), 1--67.
\bibitem{BK2} Y. Bugeaud, D. H. Kim, {\em A new complexity function, repetitions in Sturmian words, and irrationality exponents of Sturmian numbers}, preprint 	arXiv:1510.00279 [math.NT], 2015.
\bibitem{BK1} Y. Bugeaud, D. H. Kim, {\em On the $b$-ary expansions of $\log(1+1/a)$ and $e$}, preprint arXiv:1510.00282 [math.NT], Ann. Scuola Normale Superiore di Pisa, to appear.
\bibitem{ChudHouches}   G. V. Chudnovsky, {\em Rational and Pad\'e approximations to solutions of linear differential equations and the monodromy theory}, Les Houches (1979), Lecture Notes in Physics {\bf 126}, 136--169, Springer, 1980.
\bibitem{chud1} G. V. Chudnovsky, {\em On applications of diophantine approximations}, 
Proc. Natl. Acad. Sci. USA {\bf 81} (1984), 7261--7265.
\bibitem{chud2} D. V. Chudnovsky and G. V. Chudnovsky, {\em Applications of Pad\'e approximations to
Diophantine inequalities in values of G-functions}, Number theory (New York, 1983/84),
 Lecture Notes in Math. {\bf 1135}, Springer, Berlin, 1985, 9--51.
\bibitem{firi} S. Fischler, T. Rivoal, {\em On the values of $G$-functions}, 
Commentarii Math. Helv. {\bf 29}.2 (2014), 313--341. 
\bibitem{galoshkin} A. I. Galochkin, {\em Lower bounds of polynomials in the values of a certain class of analytic
functions}, in Russian, Mat. Sb. {\bf 95} (137) (1974), 396--417, 471.
\bibitem{ridout} D. Ridout, 
{\em Rational approximations to algebraic numbers}, Mathematika {\bf 4} (1957), 125--131.
\bibitem{riv1} T. Rivoal, {\em Convergents and irrationality measures of logarithms}, 
Rev. Mat. Iberoamericana {\bf 23}.3 (2007), 931--952.
\bibitem{riv3} T. Rivoal, {\em On the bits couting function of real numbers}, 
J. Austr. Math. Soc. {\bf 85} (2008), 95--111.
\bibitem{Shidlovski} 
A. B. Shidlovsky, {\em Transcendental Numbers}, de Gruyter Studies in Mathematics {\bf 12}, 1989.
\bibitem{siegel} C. Siegel,  {\em \"Uber einige Anwendungen 
diophantischer Approximationen}, Abhandlungen Akad. Berlin 1929, no. 1, 70 S.
\bibitem{zudilin2} W. Zudilin, {\em On rational approximations of values of a certain class of entire functions},  
Sb. Math. {\bf 186}.4 (1995), 555--590; translated from the 
russian version Mat. Sb. {\bf 186}.4 (1995), 89--124.
\bibitem{zudilin} W. Zudilin, {\em On a measure of irrationality for values of $G$-functions}, 
Izv. Math. {\bf 60}.1 (1996), 91--118; translated from the russian version 
 Izv. Ross. Akad. Nauk Ser. Mat. {\bf 60}.1 (1996), 87--114.
\end{thebibliography}
\end{document}